\documentclass[11pt,a4paper]{article}
\usepackage{amsmath,amsfonts,amssymb,mathrsfs,comment,soul}
\usepackage{xcolor}
\usepackage[english]{babel}
\usepackage{mathscinet}

\newtheorem{theorem}{\bf Theorem}[section]

\newtheorem{definition}[theorem]{\bf Definition}
\newtheorem{remark}[theorem]{\bf Remark}
\newtheorem{proposition}[theorem]{\bf Proposition}

\newenvironment{proof}{\noindent {\sc Proof}}{\hfill $\square$}

\newcommand{\R}{\mathbb{R}}

\def \I {{\mathbb{I}}}
\def \L {\mathscr{L}}
\def \E {\mathscr{L}_{\mathbb{G}}}
\def \H {\mathcal{H}}
\def \S {\mathcal{S}}

\def \d {{\delta}}
\def \e {{\varepsilon}}

\def \r {{\varrho}}

\def \phi {{\varphi}}

\def \G {{\Gamma}}
\def \O {{\Omega}}

\def \div {{\text{\rm div}}}
\def \loc {{\text{\rm loc}}}

\def \diag {{\text{\rm diag}}}

\def\p{\partial}

\newcommand{\scp}[1]{\langle #1 \rangle}
\newcommand{\mres}{\mathbin{\vrule height 1.6ex depth 0pt width 0.13ex\vrule height 0.13ex depth 0pt width 1.3ex}}

\usepackage{mathtools}
\DeclarePairedDelimiter{\abs}{\lvert}{\rvert}

\begin{document}
\title{Mean value formulas for classical solutions 
\\
to some degenerate elliptic equations in Carnot groups}
\author{{
\sc{Diego Pallara}
\thanks{Dipartimento di Matematica e Fisica ``Ennio De Giorgi'', Universit\`{a} del Salento and INFN, Sezione di Lecce, Ex Collegio Fiorini - Via per Arnesano - Lecce (Italy). E-mail: 	diego.pallara@unisalento.it} \ 
\sc{Sergio Polidoro}
\thanks{Dipartimento di Scienze Fisiche, Informatiche e Matematiche, Universit\`{a} degli Studi di Modena e Reggio Emilia, Via Campi 213/b, 41125 Modena (Italy). E-mail: 	sergio.polidoro@unimore.it}
}}
	
\date{ }
	
\maketitle

\begin{abstract}
We prove surface and volume mean value formulas for classical solutions to uniformly elliptic equations in divergence form with H\"{o}lder continuous coefficients. The kernels appea\-ring in the integrals are supported on the level and superlevel sets of the fundamental solution relative the adjoint differential operator. We then extend the aforementioned formulas to some subelliptic operators on Carnot groups. In this case we rely on the theory of finite perimeter sets on stratified Lie groups.
\end{abstract}
	
\setcounter{equation}{0} 
\section{Introduction}\label{secIntro}

 The mean value property that characterizes harmonic functions is a fundamental tool in potential theory and provides us with simple proofs of maximum principles, Harnack inequalities, regulari\-ty properties of harmonic functions and compactness results. 

It is known that the proof of mean value formulas relies on the divergence theorem and on the fact that the Euclidean sphere is a level set of the fundamental solution of the Laplace equation. Based on this idea, this formula has been first extended in 1951 by Pini \cite{Pini1951} to the heat equation in one space variable, then by Watson \cite{Watson1973} in 1973 to the heat equation in several space variables. Mean value formulas for degenerate operators in the form of sum of squares of H\"ormander vector fields have been proved in 1993 by Citti, Garofalo and Lanconelli \cite{CittiGarofaloLanconelli}. Still in the framework of hypoelliptic operators, we recall the article by Garofalo and Lanconelli \cite{GarofaloLanconelli-1990}, who prove in 1990 mean value formulas for Kolmogorov operators. 

The general class of uniformly parabolic operators with smooth coefficents has been considered by Fabes and Garofalo \cite{FabesGarofalo} in 1987, and by Garofalo and Lanconelli in \cite{GarofaloLanconelli-1989} in 1989. More recently, uniformly parabolic operators with weaker regularity assumptions on the coefficients of the operator have been considered by the authors and Malagoli \cite{MalPalPol}. We point out that the main difficulty encountered in \cite{MalPalPol} is due to the fact that the fundamental solution of the parabolic equation is not explicitly known and that, as a consequence, it is impossible to check the regularity of its level sets, which are the boundaries of the domain where the divergence theorem needs to be applied. This problem was bypassed in \cite{FabesGarofalo, GarofaloLanconelli-1989} by assuming the smoothness of the coefficients, which guarentees the applicability of Sard's theorem, and has been circumvented in \cite{MalPalPol} in two ways: one relies on a refinement of the classical divergence theorem in the classical vein, the other on the theory of \emph{sets with finite perimeter}. 

The aim of the present paper is to come back to elliptic operators, with natural regularity assumptions, and to address the problem for their degenerate counterpart on Carnot groups. By \emph{natural regularity} we mean that we consider possibily degenerate elliptic operators with $C^{1+\alpha}$ diffusion and drift coefficients, so that both $\L$ and its adjoint $\L^*$ can be written in non-divergence form with H\"older continuous coefficients. 

As we will see in the following, the mean value formula for \emph{uniformly elliptic operators} with natural conditions on the coefficients is rather simple. Indeed, several results on the existence of a fundamental solution $\Gamma$ of uniformly elliptic equations with H\"older continuous coefficients are avaliable in literature. All of them rely on the Levi's parametrix method. We mainly refer on the Kalf article \cite{Kalf}, where some bounds of the gradient of $\Gamma$ ensure that, for sufficiently large $C$, the level set $\big\{\Gamma = C\big\}$ is a $C^1$ manifold. We give the proof of the mean value formulas for uniformly elliptic operators  for the sake of completness, and also because it describes in the Euclidean setting the procedure that we adopt in the more difficult setting of the Carnot groups. 

The generalization of the mean value property to \emph{subelliptic operators} in a Carnot group ${\mathbb{G}}$ involves a substantial difficulty. Indeed, the fundamental solution $\Gamma$ of a subelliptic equation $\L_{\mathbb{G}} u = 0$ is not smooth, as only first and second order  \emph{horizontal} derivatives are defined as continuous functions. Therefore the level sets of $\Gamma$ may have wide subsets of singular points and the same assertion is true for the fundamental solution $\Gamma^*$ of the adjoint equation $\L^*_{\mathbb{G}} v = 0$. In this setting we cannot rely on any refinement of the classical divergence theorem as we did in \cite{MalPalPol}, and we rely on the theory of \emph{sets with finite perimeter}, which has been developed in the framework of stratified groups. We refer to the lecture notes \cite{SerCassEMS} by Serra Cassano and the wide bibliography therein, which can be applied also to families of \emph{degenerate parabolic subelliptic operators}, see \cite{PP}. 

We present our results on uniformly elliptic operators in Section \ref{SecUnif}, where we also recall the main results we
need on the existence and the main properties of a fundamental solution. We devote Section \ref{SecCarnot} to recalling the setting of Carnot groups, the relevant notions on sets of finite perimeter, an existence result on the fundamental solution, and to proving the mean value formulas for a class of subelliptic operators in Carnot groups. 

\setcounter{equation}{0} 
\section{Uniformly elliptic operators}\label{SecUnif}
Let $\Omega$ be an open subset of $\R^{N}$. We consider classical solutions $u$ to the equation $\L u = f$ in $\Omega$, where $\L$ is an elliptic operator in divergence form defined for $x \in \R^{N}$ as follows
\begin{equation} \label{e-L}
	\L u (x)  := \sum_{i,j=1}^{N}\tfrac{\p}{\p x_i} \left( a_{ij} (x) \tfrac{\p u}{\p x_j}(x) \right) + 
	\sum_{i=1}^{N} b_{i} (x) \tfrac{\p u}{\p x_i}(x) + c(x) u(x).
\end{equation}
In the following we use the notation 
$A(x) := \left( a_{ij}(x) \right)_{i,j=1,\dots, N}, b(x) := \left( b_{1} (x), \dots, b_{N} (x) \right)$ and we write $\L u $ in the short form
\begin{equation} \label{e-LL}
	\L u (x)  := \div \left( A (x) \nabla u(x) \right) + \langle b(x), \nabla u(x)\rangle + c(x) u(x).
\end{equation}
The adjoint operator is of course
\begin{equation*}
	\L^* u (x)  := \div \left( A (x) \nabla u(x) \right) - \langle b(x), \nabla u(x)\rangle + (c(x)-{\rm div}\,b(x)) u(x).
\end{equation*}
Here $\div, \nabla$ and $\langle \, \cdot \, , \, \cdot \, \rangle$ denote the divergence, the gradient and the inner product in $\R^N$, respectively. We assume that the matrix $A(x)$ is symmetric and that $\L$ is uniformly elliptic, \emph{i.e.}, there exist two constants $\lambda, \Lambda$, with $0 < \lambda < \Lambda$, such that
\begin{equation} \label{e-up}
\lambda |\xi|^2  \le \langle A(x) \xi, \xi \rangle \le \Lambda |\xi|^2, 
\end{equation}
for every $\xi \in \R^N$, for every $x \in \R^{N}$, and for $i,j=1, \dots, N$. We finally assume that the coefficients $a_{ij}, b_i, c$ and their derivatives $\frac{\p a_{ij}}{\p x_i}, \frac{\p b_{i}}{\p x_i}$ are bounded H\"older continuous functions of exponent $\alpha \in ]0,1]$, for $i,j=1, \dots, N$. Under the above assumptions, the classical parametrix method provides us with the existence of a fundamental solution $\Gamma$. Let us quote from the article of Kalf \cite{Kalf} the results we need for our purposes. We denote by $D_\Omega$ the diagonal $D_\Omega=\{(x,x)\in \Omega\times\Omega\}$.

According to the Definition at page 259 of \cite{Kalf}, we say that a function $\Gamma: (\Omega \times \Omega) \backslash D_\Omega \to \R$ is a fundamental solution to the equation $\L u = 0$ if the following conditions hold for every $y \in \Omega$:
\begin{description}
 \item[{\it i)}] $\Gamma (\cdot, y) \in C^2(\Omega \backslash \{y\})$ and $\L \, \Gamma (\cdot, y)=0$;
 \item[{\it ii)}] 
 \begin{equation*}
    \int_\Omega |\Gamma (x, y)| \, dx < + \infty;
 \end{equation*}
 \item[{\it iii)}] For every $\varphi \in C_c^\infty(\Omega)$ we have
  \begin{equation*}
    \varphi(y) = - \int_\Omega \Gamma (x, y) \L^* \varphi(x) \, dx.
 \end{equation*}
\end{description}
In the following we use property {\it iii)} with $\varphi \in C_c^2(\Omega)$, by a density argument. 

Under our hypotheses, Theorem 5 of \cite{Kalf} ensures that for every connected open set $\Omega_1$ such that $\overline{\Omega_1} \subset \Omega$ there exists a fundamental solution $\Gamma :(\Omega_1 \times \Omega_1) \backslash D_\Omega \to \R$ of $\L u = 0$. Our assumptions on the coefficients of $\L$ ensures that there exists as well a fundamental solution $\Gamma^* :(\Omega_1 \times \Omega_1) \backslash D_\Omega \to \R$ of $\L^* u = 0$.

The parametrix method used in \cite{Kalf} also provides us with the following estimates: there exist positive constants $c^-, c^+, c_0, c_1$, only depending on the dimension $N$, on the operator $\L$ and on the set $\Omega_1$, such that the following bounds hold for every $x,y \in \Omega_1$:
\begin{equation} \label{eq-bd1}
 c^- \log\left( \tfrac{\sqrt{\lambda}}{|x-y|} \right) - c_0 |x-y|^{\alpha/2} \le 
 \Gamma(x,y) \le c^+ \log\left( \tfrac{\sqrt{\Lambda}}{|x-y|} \right) + c_0 |x-y|^{\alpha/2}
\end{equation}
in the case $N=2$, while for $N>2$ we have
\begin{equation} \label{eq-bd2}
 \frac{c^-}{|x-y|^{N-2}} - \frac{c_0}{|x-y|^{N-2-\alpha}} \le \Gamma(x,y) 
 \le \frac{c^+}{|x-y|^{N-2}} + \frac{c_0}{|x-y|^{N-2-\alpha}}.
\end{equation}
Moreover, 
\begin{equation} \label{eq-bd3}
 \frac{c^-}{|x-y|^{N-1}} - \frac{c_0}{|x-y|^{N-1-\alpha}} 
 \le | \nabla_x \Gamma(x,y)| \le 
 \frac{c^+}{|x-y|^{N-1}} + \frac{c_0}{|x-y|^{N-1-\alpha}}
\end{equation}
for every $x,y \in \Omega_1$. The above inequalities also hold for $\Gamma^*$. We refer to Theorem 3 and to the formulas (2.5), (2.6) and (2.8) of \cite{Kalf}. 

We introduce some notation and we state our main results. For every $x_0 \in \R^{N}$ and for every $r>0$, we set
\begin{equation} \label{e-Psi}
\begin{split}
	\psi_r(x_0) & := \left\{ x \in \R^{N} \mid \Gamma^*(x, x_0 ) = \tfrac{1}{r^{N-2}} \right\}, 
	\\
	\Omega_r(x_0) & := \left\{ x \in \R^{N} \mid \Gamma^*(x, x_0) > \tfrac{1}{r^{N-2}} \right\}.
\end{split}
\end{equation}
Note that, as observed in Remark 7 of \cite{Kalf}, the bounds \eqref{eq-bd1} and \eqref{eq-bd2} imply $\Gamma(x,y)>0$ and $\Gamma^*(x,y)>0$ whenever $x$ is sufficiently close to $y$. Moreover, $\Omega_r(x_0)$ is bounded and $\psi_r(x_0)$ is a $C^1$ manifold, for sufficiently small positive $r$. 

We finally introduce the following kernels 
\begin{equation} \label{e-kernels}
\begin{split}
	K (x_0, x) & := 
	\frac{\langle A(x)\nabla_x\Gamma^*(x,x_0),\nabla_x\Gamma^*(x, x_0) \rangle }{|\nabla_x\Gamma^*(x,x_0)|},
	\\
	M(x_0,x) & := \frac{N}{(N-2)} \cdot
\frac{\langle A(x) \nabla_x\Gamma^*(x,x_0),\nabla_x\Gamma^*(x,x_0) \rangle}{\Gamma^*(x,x_0)^{\frac{2(N-1)}{N-2}}}.
\end{split}
\end{equation}
Note that $\nabla_x \Gamma^*(x_0,x) \neq 0$ in $\Omega_r(x_0)$ for sufficiently small $r$ and $\Gamma^*(x,x_0) \neq 0$ in $\psi_r(x_0)$, by its very definition. In the following, $\H_e^{N-1}$ denotes the $(N-1)$-dimensional Hausdorff measure, see e.g. \cite{Morgan}. The first achievements of this note are the following mean value formulas.

\begin{theorem} \label{th-1}
 Let $\Omega$ be an open subset of $\R^{N}$, with $N>2$, $f\in C(\Omega)$ and let $u$ be a classical solution to $\L u = f$ in $\Omega$. Then, for every $x_0 \in \Omega$ there is $r_0>0$ such that for every $0<r<r_0$ we have 
 \begin{equation} \label{e-meanvalue}
 \begin{split}
	u(x_0) = \int_{\psi_r(x_0)} K (x_0,x) u(x) \, d \H^{N-1} (x) 
	+ & \int_{\Omega_r(x_0)} f (x) \left( \tfrac{1}{r^{N-2}} - \Gamma^*(x,x_0) \right)\ dx  \\
	 + & \frac{1}{r^{N-2}} \int_{\Omega_r(x_0)} \left( \div \, b(x) - c(x) \right) u(x) \ dx, 
\end{split}
\end{equation}	 
\begin{align} \nonumber
	u(x_0) = \frac{1}{r^{N}} \int_{\Omega_r(x_0)} \!\!\!\!\!  M (x_0,x) &u(x) \, dx 
	+ \frac{{N}}{r^{N}}  \int_0^{r} \left(\r^{N-1}\!\!
	\int_{\Omega_\r (x_0)} \!\!\!\!\!
	f (x) \left( \tfrac{1}{\r^{N-2}} - \Gamma^*(x,x_0) \right) dx \right) d \r  \\ \label{e-meanvalueOmega}
    &+   \frac{{N}}{r^{N}} \int_0^{r} \left(\r
	\int_{\Omega_\r (x_0)} \left( \div \, b(x) - c(x) \right) u(x) \, dx \right) d \r.
\end{align}
\end{theorem}

\begin{remark}
 {\rm The usual proof of the mean value formula relies on the application of the second Green identity  on the set $\Omega_r(x_0) \backslash \Omega_\varepsilon (x_0)$, with $0 < \varepsilon < r \le r_0$, which produces the first integral in the right hand side of \eqref{e-meanvalue} and $\int_{\psi_\varepsilon(x_0)} K (x_0,x) u(x) \, d \H^{N-1} (x)$. The conclusion then follows from the fact that $\int_{\psi_\varepsilon(x_0)} K (x_0,x) u(x) \, d \H^{N-1} (x) \to u(x_0)$ as $\varepsilon \to 0.$ We adopt here an alternative approach which relies on the properties of the fundamental solution. This approach simplifies the proof of the mean value formula in the setting of degenerate subelliptic operators on Carnot groups.}
\end{remark}

\begin{proof}. Let $r_1>0$ be such that $\overline{\Omega_{r_1}(x_0)}\subset\Omega$ and consider an open set $\Omega_1$ such that $\overline{\Omega_{r_1}(x_0)} \subset \Omega_1$ and that $\overline{\Omega_1} \subset \Omega$. Then there exists the fundamental solution $\Gamma^*$ of the equation $\L^* v = 0$ in $\Omega_1$. As noticed above, we can, and we do, choose $r_0 \in ]0, r_1]$ such that $|\nabla_x\Gamma^*(\cdot,x_0)|\neq 0$ in $\Omega_{r_0}(x_0)$, hence $\psi_r(x_0)$ is a $C^1$-manifold for every positive $r < r_0$. We choose $r \le r_0$ and we introduce a further positive parameter $\varepsilon$ small enough to have $B_\varepsilon(x_0):=\{|x-x_0|<\varepsilon\}\subset\Omega_r(x_0)$. We let $\phi_\varepsilon \in C^\infty(\R^N)$ be a function such that ${\rm supp}(\phi_\varepsilon) \subset B_\varepsilon(x_0)$ and $\phi_\varepsilon(x) = 1$ for $|x-x_0|\leq \varepsilon/2$. Since $\phi_\varepsilon u \in C_c^2(\Omega_1)$, the very definition of $\Gamma^*$ yields
\begin{equation*}
\phi_\varepsilon (y) u (y) = - \int_\Omega \Gamma^* (x,y) \L (\phi_\varepsilon u) (x) \, dx
\end{equation*}
for every $y \in \Omega_1$. In particular,
\begin{equation} \label{eq-mvf-1}
     u (x_0) = - \int_\Omega \Gamma^* (x,x_0) \L (\phi_\varepsilon u) (x) \, dx.
\end{equation}
We next consider the functions $w(x) := (1 - \phi_\varepsilon(x) )u(x)$ and $v(x) := \Gamma^*(x,x_0) - \frac{1}{r^{N-2}}$ and we note that 
\begin{equation} \label{eq-div-L*}
\begin{split}
 w(x) \L^* v(x) - v(x) \L w(x) = 
 & \div \big( w(x) A(x) \nabla v(x) - v(x) A(x) \nabla w(x) \big)  \\
 & - \div \big( w(x) v(x) b(x) \big)
\end{split}
\end{equation}
for every $x \in \Omega \backslash \big\{x_0 \big\}$. Recalling that $\L^* v = \frac{1}{r^{N-2}} \left( \div \, b - c \right)$ in 
$\Omega \backslash \big\{x_0 \big\}$, equation \eqref{eq-div-L*} can be written as follows
\begin{equation*}
 \frac{1}{r^{N-2}} \left( \div \, b(x) - c(x) \right) w(x) - v(x) \L w(x) = 
 \div \, \Phi (x), \quad \Phi(x) := \big( w A \nabla v - v A \nabla w - wv b\big)(x).
\end{equation*}
By our choice of $\phi_\varepsilon$ we have that $\Phi$ vanishes in 
$B_{\varepsilon/2}(x_0)$, then it can be extended to a $C^1(\Omega)$ function by setting $\Phi(x_0)= 0$. By the divergence theorem we find
\begin{equation} \label{eq-div-k}
 \int_{\Omega_r(x_0)} \left( \tfrac{1}{r^{N-2}} \left( \div \, b(x) - c(x) \right) w(x) - v(x) \L w(x)\right) dx
=  - \int_{\psi_r(x_0)} \scp{\nu,\Phi} d \H_e^{N-1},
\end{equation}
where $\nu(x)=\tfrac{\nabla_x \Gamma^*(x,x_0)}{\abs{\nabla_x \Gamma^*(x,x_0)}}$. 
In order to conclude our proof, we rewrite \eqref{eq-mvf-1} in its equivalent form
\begin{equation} \label{eq-div-k2}
 u (x_0) = - \int_{\Omega_r(x_0)} \left( \tfrac{1}{r^{N-2}} + v(x) \right) \L (\phi_\varepsilon u) (x) dx.
\end{equation}
We then recall that $u = w + \phi_\varepsilon u$, so that $f = \L u = \L w + \L (\phi_\varepsilon u)$. Than, by subtracting \eqref{eq-div-k2} and \eqref{eq-div-k}, we find
\begin{equation} \label{eq-div-eps}
\begin{split}
u(x_0) =& - \int_{\Omega_r(x_0)} v(x)f(x)\, dx
+ \int_{\Omega_r(x_0)} \left( \tfrac{1}{r^{N-2}} ( \div \, b(x) - c(x)) \right) (1-\phi_\varepsilon(x))u(x)\, dx
\\
& + \int_{\psi_r(x_0)} \scp{\nu,\Phi} d \H_e^{N-1}
- \tfrac{1}{r^{N-2}} \int_{B_\varepsilon(x_0)} \L(\phi_\varepsilon u)(x)\, dx.
\end{split}
\end{equation}
Now we let $\varepsilon\to 0$. Then, $\phi_\varepsilon(x)\to 0$ a.e. and moreover
\[
\int_{B_\varepsilon(x_0)} \L(\phi_\varepsilon u)(x)\, dx
= \int_{B_\varepsilon (x_0)} {\rm div}\, (A\nabla(\phi_\varepsilon u)) \, dx
+\int_{B_\varepsilon (x_0)} \langle b,\nabla(\phi_\varepsilon u)\rangle \, dx 
+ \int_{B_\varepsilon (x_0)} c\phi_\varepsilon u\, dx.
\]
The first integral vanishes by the divergence theorem, as $\phi_\varepsilon=0$ on the boundary, the last one tends to $0$ because the function $\phi_\varepsilon u$ is bounded, whereas for the second one we notice that
\[
\int_{B_\varepsilon (x_0)} \langle b,\nabla(\phi_\varepsilon u)\rangle \, dx 
= \int_{\partial B_\varepsilon (x_0)} \phi_\varepsilon u \langle b,\nu\rangle \, d\H_e^{N-1}
- \int_{B_\varepsilon (x_0)} \phi_\varepsilon u {\rm div}\, b \, dx
\]
where the surface integral vanishes and the second one tends to $0$. Finally, we have  
\begin{equation*}
 \scp{\nu,\Phi}(x) = K(x_0,x)u(x),
\end{equation*}
since $\nu(x)=\tfrac{\nabla_x \Gamma^*(x,x_0)}{\abs{\nabla_x \Gamma^*(x,x_0)}}$. This concludes the proof of the first statement of Theorem \ref{th-1}. 

The proof of the second assertion of Theorem \ref{th-1} is a direct consequence of the first one and of the coarea formula. Indeed, fix a positive $r$ as above, multiply \eqref{e-meanvalue} (with $\r$ in place of $r$) by $\frac{N}{r^{N}}\r^{N-1}$ and integrate over $]0,r[$. We find
\begin{equation*} 
\begin{split}
	\frac{N}{r^{N}} \int_0^r \varrho^{N-1} u(x_0) d \varrho = & 
	\frac{N}{r^{N}} \int_0^r \varrho^{N-1} \bigg(\int_{\psi_\varrho(x_0)} 
	K (x_0,x) u(x) \, d \H_e^{N-1} (x) \bigg) d \varrho\, 
	\\
	& + \frac{N}{r^{N}} \int_0^r \varrho^{N-1} \bigg( \int_{\Omega_\r(x_0)} f (x) 
	\left( \tfrac{1}{\r^{N-2}} - \Gamma^*(x,x_0) \right) dx \bigg) d \varrho  
	\\
	& + \frac{N}{r^{N}} \int_0^r \varrho \bigg( 
	\int_{\Omega_\r (x_0)} \left( \div \, b(x) - c(x) \right) u(x) \, dx \bigg) d \varrho.
\end{split}
\end{equation*}
The left hand side of the above equality equals $u(x_0)$, while the last two terms agree with the last two terms appearing in the statement of Theorem \ref{th-1}. In order to conclude the proof we only need to show that 
\begin{equation} \label{e-meanvalue-step2}
\frac{1}{N-2} \int_{\Omega_r (x_0)} M (x_0,x) u(x)\, dx =
\int_0^r \varrho^{N-1} \bigg(\int_{\left\{\Gamma^*(\cdot,x_0)= \tfrac{1}{\varrho^{N-2}}\right\} } 
K (x_0,x) u(x) \, d \H_e^{N-1} (x) \bigg) d \varrho.
\end{equation}
With this aim, we substitute $y = \frac{1}{\varrho^{N-2}}$ in the left hand side of \eqref{e-meanvalue-step2} and we recall the definition of the kernel $K$. We find 
\begin{align} \label{e-meanvalue-step3}
	  \int_0^r \varrho^{N-1} & 
\bigg( \int_{\left\{\Gamma^*(\cdot,x_0)= \tfrac{1}{\varrho^{N-2}}\right\} } 
	  \frac{\langle A(x)\nabla_x\Gamma^*(x,x_0),\nabla_x\Gamma^*(x, x_0) \rangle }{|\nabla_x\Gamma^*(x,x_0)|}
	  u(x) \, d \H^{N-1}_e (x) \bigg) d \varrho 
\\ \nonumber
     & \!\!\!\!\!\!\!\!\!\! = \frac{1}{N-2}
   \int_{\frac{1}{r^{N-2}}}^{\infty} \frac{1}{y^{\frac{2(N-1)}{N-2}}} \bigg(\int_{\left\{\Gamma^*(\cdot,x_0) = y\right\} } 
	\frac{\langle A(x)\nabla_x\Gamma^*(x,x_0),\nabla_x\Gamma^*(x, x_0) \rangle }{|\nabla_x\Gamma^*(x,x_0)|}
	u(x) \, d \H_e^{N-1} (x) \bigg) d y 
	\\ \nonumber
	& \!\!\!\!\!\!\!\!\!\! =\frac{1}{N-2} \int_{\frac{1}{r^{N-2}}}^{\infty} \bigg(\int_{\left\{\Gamma^*(\cdot,x_0) = y\right\} }\!\!\!
	\frac{\langle A(x)\nabla_x\Gamma^*(x,x_0),\nabla_x\Gamma^*(x, x_0) \rangle }
	{\Gamma^*(x,x_0)^{\frac{2(N-1)}{N-2}}|\nabla\Gamma^*(x,x_0)|} u(x) \, d \H^{N-1}_e (x) \bigg) d y.
\end{align}
We conclude the proof of \eqref{e-meanvalue-step2} by applying the coarea formula, see e.g. \cite{Morgan}.
\end{proof}

\begin{remark}\label{rm-2}{\rm 
If $N=2$ we start of course from \eqref{eq-bd1}. Let us show how the results and the proofs have to be adapted. First, the integration domains are
\begin{equation} \label{e-Psi2}
\begin{split}
	\psi_r(x_0) & := \left\{ x \in \R^{2} \mid \Gamma^*(x, x_0 ) = \log\left(\tfrac{1}{r}\right) \right\}, 
	\\
	\Omega_r(x_0) & := \left\{ x \in \R^{2} \mid \Gamma^*(x, x_0) > \log\left(\tfrac{1}{r}\right) \right\},
\end{split}
\end{equation}
the kernel $K(x_0,x)$ is the same as in the case $N\geq 3$, whereas
\begin{equation} \label{e-kernels2}
M(x_0,x) := 2
\frac{\langle A(x) \nabla_x\Gamma^*(x,x_0),\nabla_x\Gamma^*(x,x_0) \rangle}{\exp\{2\Gamma^*(x,x_0)\}}.
\end{equation}
Arguing as in the case $N\geq 3$ we easily get 
 \begin{equation} \label{e-meanvalue2}
 \begin{split}
	u(x_0) = \int_{\psi_r(x_0)} K (x_0,x) u(x) \, d \H^{1} (x) 
	+ & \int_{\Omega_r(x_0)} f (x) \left( \log\left( \tfrac{1}{r}\right) - \Gamma^*(x,x_0) \right)\ dx  \\
	 + & \log\left(\frac{1}{r}\right) 
\int_{\Omega_r(x_0)} \left( \div \, b(x) - c(x) \right) u(x) \ dx\,d\r. 
\end{split}
\end{equation}	 
In order to deduce the mean value formula corresponding to 
\eqref{e-meanvalueOmega} we write
\begin{align*}
u(x_0) =& \frac{2}{r^2}\int_0^r \r u(x_0)\,d\r = 
\frac{2}{r^2} \int_0^r \r \int_{\psi_\r(x_0)} 
K (x_0,x) u(x) \, d \H^{1} (x)\, d\r 
\\
&+\frac{2}{r^2}\int_0^r \r \int_{\Omega_\r(x_0)} f (x) 
\left( \log\left( \tfrac{1}{r}\right) - \Gamma^*(x,x_0) \right)\ dx\, d\r  
\\
&+\frac{2}{r^2}\int_0^r \r \log\left(\frac{1}{\r}\right) 
\int_{\Omega_\r(x_0)} \left( \div \, b(x) - c(x) \right) u(x) \ dx\, d\r.
\end{align*}	 
Concerning the first integral, substituting $y=\log(1/\r)$ we have 
\begin{align*}
&\int_0^r \r \int_{\{\Gamma^*(\cdot,x_0)=\log(1/\r\}} 
\frac{\langle A(x)\nabla_x\Gamma^*(x,x_0),\nabla_x\Gamma^*(x, x_0) \rangle }{|\nabla_x\Gamma^*(x,x_0)|} u(x) \, d \H^{1} (x)\, d\r 
\\
&\int_{\log(1/r)}^{\infty} e^{-2y} \int_{\{\Gamma^*(\cdot,x_0)=y\}} 
\frac{\langle A(x)\nabla_x\Gamma^*(x,x_0),\nabla_x\Gamma^*(x, x_0) \rangle }{|\nabla_x\Gamma^*(x,x_0)|} u(x) \, d \H^{1} (x)\, dy
\\ 
& \int_{\Omega_r(x_0)} \frac{\langle A(x)\nabla_x\Gamma^*(x,x_0),\nabla_x\Gamma^*(x, x_0) \rangle }{\exp\{2\Gamma^*(x,x_0)\}} u(x) \, dx,
\end{align*}
whence
\begin{align*}
u(x_0)=& \frac{1}{r^2}\int_{\Omega_r(x_0)} M(x_0,x) u(x)\, dx
\\
&+\frac{2}{r^2}\int_0^r \r \int_{\Omega_\r(x_0)} f (x) 
\left( \log\left( \tfrac{1}{r}\right) - \Gamma^*(x,x_0) \right)\ dx\, d\r  
\\
&+\frac{2}{r^2}\int_0^r \r \log\left(\frac{1}{\r}\right) 
\int_{\Omega_\r(x_0)} \left( \div \, b(x) - c(x) \right) u(x) \ dx.
\end{align*}	 
}\end{remark}

\begin{remark}
{\rm In Theorem \ref{th-1} and in Remark \ref{rm-2} we have assumed $r_0$ small enough in order to exploit the regularity of the level sets of $\Gamma^*$ and apply the classical divergence theorem. Indeed, as we see in the next section, the mean value formulas hold true for almost every $r$ such that $\Omega_r(x_0)$ is a bounded open subset of $\Omega$. This could be proved, in the same vein, relying on the theory of sets with finite perimeter in $\R^N$.}
\end{remark}

\setcounter{equation}{0} 
\section{Subelliptic operators in Carnot groups}\label{SecCarnot}
In this section we state and prove the mean value formula for a class of subellitic operators in Carnot groups. We need much more preliminar information with respect to the uniformly elliptic case treated in Section \ref{SecUnif}, and accordingly we split this section in various subsections. In the first one we describe the structure of the Carnot groups. In the second one we present the class of operators we are interested in, we recall the properties of fundamental solutions and state the main result, Theorem \ref{th-2}. In the third one we discuss the properties of sets with finite perimeter and, finally, in the last one we prove Theorem \ref{th-2}.

We consider $m$ smooth vector fields ($1\leq m\leq N$) 
\begin{equation}\label{defX}
     X_{j}(x)=\sum_{k=1}^{N} \varphi_{k}^{j}(x)\p_{x_{k}}, \qquad j=1,\dots,m,
\end{equation}
with $\varphi_k^j\in C^\infty(\R^N)$. We introduce the Lie algebra generated by $X_1, \ldots, X_m$  
\begin{equation} \label{e-Lie}
    {\mathfrak{g}} = {\rm Lie} (X_1, \ldots, X_m)
\end{equation}
and we assume the following:
\begin{description}
\item[{\rm [H.1]}] The vector fields $X_1,\ldots,X_m$ satisfy the H\"{o}rmander's rank condition
\begin{equation} \label{e-hrc}
    {\rm rank} \, {\mathfrak{g}} (x) = N \qquad \text{for every} \quad x \in \R^{N}.
\end{equation}
  \item[{\rm [H.2]}] there exists a homogeneous Lie group
  $\mathbb{G}=\left(\R^N,\circ, \d_{\lambda}\right)$ such that
\begin{description}
  \item[{\it i)}] $X_{1},\dots,X_{m}$ are left translation invariant on   $\mathbb{G}$;
  \item[{\it ii)}] $X_{1},\dots,X_{m}$ are $\d_{\lambda}$-homogeneous of degree one.
\end{description}
Moreover, $X_j(0)$ agrees with the $j$-th element of the canonical basis of $\R^N$, for $j=1, \dots, m$. 
\end{description}

\subsection{Stratified groups}\label{subsGroups}
A Lie group $\mathbb{G}=\left({\mathbb R}^N,\circ\right)$ is said {\it homogeneous} if a family of dilations $\left(\d_{\lambda}\right)_{\lambda>0}$ exists on $\mathbb{G}$ and it is an automorphism of the group:
\begin{equation*}
\d_\lambda(x \circ y) = \left(\d_\lambda x\right) \circ \left( \d_\lambda y\right), \quad \text{for all} \ x,y
\in \R^{N} \ \text{and} \ \lambda >0.
\end{equation*}
The assumptions [H.1] and [H.2] induce a direct sum decomposition of the Lie algebra $\mathfrak{g}$
\begin{equation}\label{e-Oplus}
    \mathfrak{g} = V_{1}\oplus\dots\oplus V_{\nu},
\end{equation}
where $V_1 = \text{span} \big\{X_1, \dots X_m\big\},  V_{k+1} = \text{span} \big\{[X, Y], \mid X \in V_1, Y \in V_k \big\},$ for $k=1, \dots, \nu-1$ and $\big\{[X, Y], \mid X \in V_1, Y \in V_{\nu} \big\} = \big\{ 0 \big\}.$
In the sequel we denote by $n_j$ the dimension of $V_j$, for $j=1,\dots, \nu$. Note that [H.2] yields $m = n_1$. The dilation $\delta_\lambda$ on $\R^{N}$ will be represented by a diagonal matrix, which necessarily has the following form
\begin{equation} \label{e-delta}
  \delta_\lambda = \diag ( \lambda \I_{n_1} , \lambda^2 \I_{n_2}, \ldots, \lambda^{\nu} \I_{n_\nu}),
\end{equation}
the integer $Q= n_1 + \dots + n_\nu$ will be called \emph{homogeneous dimension} of $\mathbb{G}$ and we have  
\begin{equation}\label{e-Q}
    \det \delta_\lambda = \lambda^{Q}.
\end{equation}
We refer to the monograph \cite{LibroBLU} for a more detailed treatment of homogeneous Lie groups and for an exhaustive bibliography on this subject. In particular, from Proposition 1.38 in \cite{LibroBLU} it follows that
\begin{equation}\label{e-adj}
  X_j^* = - X_j, \qquad j=1, \dots, m. 
\end{equation} 
Let us introduce the distance that we use in this paper. From \cite[Theorem 5.1]{fraserser3} we know that there are constants $\varepsilon_j\in ]0,1],\ j=1,\ldots,\nu$, with $\varepsilon_1=1$, such that the function 
\begin{equation}\label{defnorma_infty}
x\mapsto \|x\|_\infty = \max_{j=1,\ldots,\nu}\{\varepsilon_j|x_j|^{1/j}\},
\end{equation}
where $x_j\in \R^{n_j}$ and $|\cdot|$ denotes the usual Euclidean norm, defines a norm and as a consequence the distance  
\begin{equation}\label{defd_infty}
d_\infty(x,y)=\|y^{-1}\circ x\|_\infty.
\end{equation}
We notice that $d_\infty$ is equivalento to the 
Carnot-Carath\'{e}odory distance and that for every compact set $K \subset \R^{N}$ there exist two positive constants $c_K^-$ and $c_K^+$, such that
\begin{equation}\label{e-dist-eucl}
  c_K^- |x-y| \le d_\infty(x,y) \le c_K^+ |x-y|^\frac{1}{\nu}, \qquad \text{for all} \ x,y \in K.
\end{equation}
The invariance properties 
\begin{equation} \label{e-dist-INV}
   d_\infty (y \circ x, y \circ z) =  d_\infty (x,z), \quad 
   d_\infty (\delta_\lambda x, \delta_\lambda y) = \lambda d_\infty (x,y),
\end{equation}
hold for every $x,y,z$ in $\R^{N}$ and for every positive $\lambda$, see again \cite{fraserser3}. For every $\nu\in V_1$, denote by $\nu^\bot$ the codimension 1 subspace of $V_1$ orthogonal to $\nu$ and introduce the hyperplane $N=\nu^\bot\oplus V_2\oplus \cdots \oplus V_\nu$ in $\R^N$ and the constant
\begin{equation}\label{deftheta}
\theta=\max_{z\in B(0,1)}\{\H_e^{N-1}(B(z,1)\cap N)\}
\end{equation}
which is independent of $\nu$ because $d_\infty$ is \emph{vertically symmetric} according to Definition 6.1 in \cite{MagnaniIUMJ}, see Remark 6.2 and Theorem 6.3 in \cite{MagnaniIUMJ}. The constant $\theta$ is introduced in \cite{MagnaniCalcVar}, is called spherical factor and is denoted $\omega_{{\mathbb G},Q-1}$ there. Let us set ${\rm diam}_{\mathbb G}(E)=\sup_{x,y\in E}d_\infty(x,y)$ and define the $h$-dimensional spherical Hausdorff measure $\S^{h}_{\mathbb G}$ of a Borel set $E$, $0\leq h\leq Q$. First set for $r>0$
\[
{\mathcal S}_{{\mathbb G},r}^h(E)= 
\inf\left\{
\sum_{i=0}^\infty \frac{\theta}{2^{h}}{\rm diam}_{\mathbb G}\,(B_i)^h: B_i\  d_\infty{\rm -balls,}\ E\subset \bigcup_{i=0}^\infty B_i,
\ {\rm diam}_{\mathbb G}\,(B_i)\leq r \right\}
\] 
and then 
\begin{equation}\label{HausdorffMeas}
{\mathcal S}_{\mathbb G}^h(E)=\lim_{r\downarrow 0}
{\mathcal S}_{{\mathbb G},r}^h(E)=\sup_{r\geq 0}
{\mathcal S}_{{\mathbb G},r}^h(E).
\end{equation}
The role of the constant $\theta$ in the definition of $\S^{h}_{\mathbb G}$ will be discussed in the Remark \ref{thetacostante} below.

Finally, we recall the notations of Lie derivatives and H\"older spaces. For any $x_0 \in \O$ and $j=1, \dots, m$, let $\gamma$ be the solution to the Cauchy problem 
\begin{equation*}
 \gamma'(s) = X_j(\gamma(s)), \qquad \gamma(0)= x_0.
\end{equation*}
Then the Lie derivative $X_j u(x_0)$ of $u$ at $x_0$ is
\begin{equation*}
 X_j u(x_0) := \frac{d}{d s}u(\gamma(s))_{\mid s = 0}.
\end{equation*}
We say that $u$ belongs to $C^{1}_{\mathbb G}(\O)$ if $u, X_j u$ for $j= 1, \dots, m$ are continuous functions on $\Omega$ and that $u\in C^{2}_{\mathbb G}(\O)$ if $u, X_j u, X_i X_j u$ are continuous in $\Omega$ for $i,j=1,\ldots,m$. For $\alpha \in ]0,1]$ we say that $u$ is $\alpha$--H\"older continuous, and we write $u \in C_{\mathbb G}^\alpha(\Omega)$, if there exists a positive constant $M$ such that
\begin{equation} \label{yeq-C-alpha}
|u(x)-u(y)| \le M d_\infty(x,y)^\alpha, \quad \text{for every} \ x,y \in \Omega.
\end{equation}
Moreover, we say that $u$ belongs to $C_{\mathbb G}^{1+ \alpha}(\O)$ (resp. $C_{\mathbb G}^{1+ \alpha}(\O)$) if $u$, the derivatives $X_1 u, \dots,  X_{m} u$ (resp. and $X_i X_j u$, $i,j= 1, \dots, m$) belong to the space $C_{\mathbb G}^{\alpha}(\O)$. The function $u$ belongs to $C^\alpha_{{\mathbb G},\loc}(\Omega)$ ($C^{k+ \alpha}_{{\mathbb G},\loc}(\O)$, respectively) if it belongs to $C_{\mathbb G}^\alpha(K)$ (resp. $C_{\mathbb G}^{k+ \alpha}(K)$) for every compact set $K \subset \O$. 

\subsection{Mean value formulas in Carnot groups}\label{IntroSublaplacians}

Let $X_1,\ldots,X_m$ be a family of vector fields satisfying the assumptions {\rm [H.1]} and {\rm [H.2]}. We consider the class of operators given as follows
\begin{equation} \label{eq-E}
\begin{split}
\E u &= \sum_{i,j=1}^m a_{ij}X_iX_ju 
     + 2 \sum_{i,j=1}^m X_ja_{ij}X_iu
 + \sum_{i,j=1}^m X_iX_ja_{ij} u
 \\
 & = {\rm div}_{\mathbb G}(A\nabla_{\mathbb G}u) 
 + \langle b, \nabla_{\mathbb G}u\rangle + cu , 
\end{split}
\end{equation}
where we have set 
\[
 b_i = \sum_{j=1}^m X_j a_{ij}, \qquad 
 c = \sum_{i,j=1}^m X_iX_j a_{ij}.
\]
Concerning the horizontal gradient $\nabla_{\mathbb G}$ and the divergence ${\rm div}_{\mathbb G}$ appearing in the above formulas, we agree to identify a \emph{horizontal section} $F=\sum_{j=1}^m F_jX_j$ with its canonical coodinates $F= \left(F_1, \dots, F_m\right)$. With this agreement, we denote the gradient of $f\in C^1_{\mathbb G}(\R^N)$ and the divergence of $F\in C^1_{\mathbb G}(\R^N,\R^m)$ by  
\begin{equation}\label{defDivGrad}
\nabla_{\mathbb G}f: =\sum_j(X_jf)X_j \quad{\rm and}\quad  
{\rm div}_{\mathbb G}F := - \sum_{j=1}^m X_j^*F_j = \sum_{j=1}^m X_jF_j.
\end{equation}

We assume that $A = (a_{ij})_{i,j=1, \dots,m}$ is a symmetric matrix satisfying the condition \eqref{e-up} for every $\xi\in\R^m$ and $x\in\R^N$ and that for every $i, j = 1, \dots, m$, the coefficients $a_{ij}$ belong to the space $C_{\mathbb G}^{2+\alpha}(\R^N)$ for some $\alpha \in ]0,1]$. The reason for this particular choice of the coefficients $b_1, \dots, b_m$ and $c$ is that, because of  \eqref{e-adj}, the \emph{adjoint operator} of $\E$ has the following simple form
\begin{equation} \label{eq-E^*}
\E^* = \sum_{i,j=1}^m a_{ij}X_iX_j.
\end{equation}
We point out that we rely on a result by Bonfiglioli, Lanconelli and Uguzzoni \cite{BLU2}, who prove the existence of a fundamental solution $\Gamma^*$ for operators $\E^*$ in the form \eqref{eq-E^*}. Less restrictive assumptions on $b_1, \dots, b_m$ and $c$ would be allowed as soon as more general existence results for $\Gamma^*$ will be available.

Let $\Omega$ denote any open subset of $\R^{N}$, and  let $u$ be a real valued function defined on $\Omega$. We say that $u$ is a \emph{classical solution} to the equation $\E u = f$ in $\Omega$ if $u$ belongs to $C_{\mathbb G}^{2}(\O)$ and the equation $\E u = f$ is satisfied at every point of $\O$. The analogous meaning is given to $\E^* v = g$.

Let us collect the relevant results on the fundamental solutions. If $\Omega$ denotes any open subset of $\R^{N}$, we let 
$D_\Omega = \{(x,x) \in \Omega \times \Omega\}$ be its diagonal. As we already pointed out in the Introduction, we assume the existence of a local fundamental solution $\Gamma^*$ to the adjoint equation $\E^* v = 0$. With this we mean a function $\G^* = \G^*(y, x)$ defined in $\big(\R^N \times \R^N \big) \backslash D_{\R^N}$, which satisfies the following conditions:
\begin{enumerate}
\item For every $x \in \R^N$ the function $\G^*( \, \cdot \, , x)$ belongs to $C_{\mathbb G}^{2}(\R^N \backslash \{ x \})$ and is a classical solution to $\E^* \, \G^*(\cdot , x) = 0$ in $\R^N \backslash \{ x \}$;
\item for every $\phi \in C_c^\infty(\O)$ the function
\begin{equation}\label{Def-Gamma}
  w(y)=\int_{\R^N}\Gamma^*(y,x)\phi(x)d x 
\end{equation}
belongs to $C_{{\mathbb G},{\rm loc}}^{2}(\R^N)$ and is a classical solution to $\E^* w = - \varphi$ in $\R^N$. 
\end{enumerate}

The existence of a fundamental solution for the operator $\E^*$ when the coefficients are constant has been proved by Folland \cite{Folland} and by Kogoj and Lanconelli in \cite{KogojLanconelli2} with a different approach. The existence of a local fundamental solution for operators $\E^*$ with H\"older continuous coefficients has been established by Bonfiglioli, Lanconelli and Uguzzoni in \cite{BLU2} for the operator $\E^*$ in non-divergence form \eqref{eq-E^*}. Let us give the complete statement of the result in \cite{BLU2}, that has been obtained by the Levi's parametrix method. 

\begin{theorem} \label{th-Gamma} {\rm (Theorem 1.5 in \cite{BLU2})} 
Let $X_1, \dots X_m$ be a family of H\"ormander vector fields satisfying the assumptions {\rm [H.1]} and {\rm [H.2]} in $\R^N$ with $N>2$. Consider the differential operator $\E^*$ in \eqref{eq-E^*}, where $A = (a_{ij})_{i,j=1, \dots,m}$ is a symmetric matrix satisfying the condition \eqref{e-up} for every $\xi\in\R^m$ and $x\in\R^N$ and for some constants $0 < \lambda < \Lambda$. Suppose that for every $i, j = 1, \dots, m$, the coefficients $a_{ij}$ belong to the space $C_{\mathbb G}^{\alpha}(\R^N)$ for some $\alpha \in ]0,1]$. Then, for every bounded open set $\Omega \subset \R^N$ there exists a fundamental solution $\G^*$ of $\E^*$ in $\Omega$. Moreover, for every $x \in \O$ the function $\G^*( \, \cdot \, , x)$ belongs to $C_{\mathbb G}^{2+\alpha}(\O \backslash \{ x \})$ and, for every compact set $K \subset \O$ there exists a positive constant $C$ such that 
\begin{equation}\label{eq-E-bds}
 0 \le \G^*(x,y) \le C \left( 1 + d_\infty(x,y)^{2-Q} \right), 
\end{equation}
for every $y \in K$ and $x \in \O$.
\end{theorem}

\begin{remark} \label{rem-L^*}{\rm Notice that the case $N \le 2$ reduces to the Euclidean one, which is not considered in this section devoted to degenerate subelliptic operators. In particular, we always have $ Q > N\geq 3$. Moreover, in the proof of our main result we use the following property of $\Gamma^*$, which follows from \eqref{Def-Gamma}. If $u \in C_{{\mathbb G},c}^{2}(\O)$ then 
\begin{equation} \label{eq-rem-L^*}
    u(y) = - \int_\Omega \Gamma^* (x,y) \E u(x) \, dx.
\end{equation}
Indeed, for $\varphi\in C_c^\infty(\Omega)$ and $u\in C^2_{{\mathbb G},c}(\Omega)$ we have 
\begin{align*}
\int_\Omega \varphi(y)\int_\Omega\Gamma^*(x,y)\E u(x)dx\ dy 
&=\int_\Omega \E u(x)\int_\Omega\Gamma^*(x,y)\varphi(y)dy\ dx 
\\
&=\int_\Omega u(x)\E^*\Bigl(\int_\Omega\Gamma^*(x,y)\varphi(y)dy\Bigr)\ dx
= - \int_\Omega u(x)\varphi(x)\, dx,
\end{align*}
and we conclude by the arbitrariness of $\varphi$. 
}\end{remark}

We are now in a position to state the main result of this section. We keep the notation $\Omega_r(x_0),\ \psi_r(x_0)$ used in Section \ref{SecUnif} and define the kernels
\begin{equation} \label{e-kernels-g}
\begin{split}
	K_{\mathbb G} (x_0, x) & := 
	\frac{\langle A(x)\nabla_{\mathbb G}\Gamma^*(x,x_0),\nabla_{\mathbb G}\Gamma^*(x,x_0) \rangle }{|\nabla\Gamma_{\mathbb G}^*(x,x_0)|},
	\\
	M_{\mathbb G}(x_0,x) & := \frac{Q}{(Q-2)} \cdot
	\frac{\langle A(x) \nabla_{\mathbb G}\Gamma^*(x,x_0),\nabla_{\mathbb G}\Gamma^*(x,x_0) \rangle}
	{\Gamma^*(x,x_0)^{\frac{2(Q-1)}{Q-2}}},
\end{split}
\end{equation}
where $\nabla_{\mathbb G}$ is defined in \eqref{defDivGrad}. We agree to set $K_{\mathbb G}(x_0,x) = 0$ whenever $\nabla_{\mathbb G} \Gamma^*(x,x_0)=0$. 

\begin{theorem} \label{th-2}
Let $\Omega$ be an open subset of $\R^{N}$, $f\in C(\Omega)$ and let $u$ be a classical solution to $\E u = f$ in $\Omega$. Then, for every $x_0 \in \Omega$ and for almost every $r>0$ such that $\overline{\Omega_r(x_0)} \subset \Omega$ we have 
 \begin{equation} \label{e-meanvalue-g}
	u(x_0) = \int_{\psi_r(x_0)} K_{\mathbb G} (x_0,x) u(x) \, 
	d\S^{Q-1}_{\mathbb G} (x) 
	+ \int_{\Omega_r(x_0)} f (x) \left( \tfrac{1}{r^{Q-2}} - \Gamma^*(x,x_0) \right)\ dx  ,
\end{equation}
\begin{equation}\label{th-2-2}
\begin{split}
	u(x_0) =& \frac{1}{r^{Q}} \int_{\Omega_r(x_0)}
	M_{\mathbb G} (x_0,x) u(x) \, dx 
	\\
	&+ \frac{{Q}}{r^{Q}} \int_0^{r} \left(\r^{Q-1}
	\int_{\Omega_\r (x_0)} f (x) \left( \tfrac{1}{\r^{Q-2}} - \Gamma(x,x_0) \right) dx \right) d \r  .
	\end{split}
\end{equation}
The second statement holds for \emph{every} $r>0$ such that ${\Omega_r(x_0)} \subset \Omega$.
\end{theorem}

\subsection{Sets of finite perimeter in stratified groups}
\label{subsPerimeters}

In this subsection we present the basic results on functions of bounded variation and sets with finite perimeter that we need to deal with fundamental solutions. If $\mu$ is a Borel measure and $E$ is a Borel set, we use the notation $\mu\mres E(B)=\mu(E\cap B)$. In the following $B(x,r)$ denotes the ball of center $x$ and radius $r$ of the distance $d_\infty$ defined in \eqref{defd_infty}. The space $BV({\mathbb{G}})$ of functions of bounded variation in $\mathbb{G}$ goes back to \cite{CapDanGar94The} and we refer to \cite{SerCassEMS} and to \cite{fraserser3} for more information. 

\begin{definition}\label{defBV} 
Let $\Omega\subset\R^N$ be an open set and $u\in L^{1}(\Omega)$; we define
\begin{equation}  \label{deftotvarG}
\left\Vert \nabla_{\mathbb G}u\right\Vert \left(\Omega\right)
=\sup\left\{  \int_{\Omega}u(x)\mathrm{div}_{\mathbb G}g(x)dx:\ g\in C_{c}^{1}\left(\Omega,{\mathbb{R}}^{m}\right),\left\Vert g\right\Vert _{\infty}\leq1\right\},
\end{equation}
where $\mathrm{div}_{\mathbb G}$ is defined in \eqref{defDivGrad}. We say that 
$u\in BV_{\mathbb{G}}(\Omega)$ if $\|\nabla_{\mathbb G}u\|(\Omega)$ is finite.
\end{definition}

\begin{remark}\label{Xdependence}{\rm
We point out (see \cite[Remarks 2.10, 2.19]{fraserser3}) that the (usual) notation $\|\nabla_{\mathbb G}u\|$ is somehow misleading, as the total variation depends upon the fixed vector fields $X_j$, even though the functional class $BV_{\mathbb G}(\Omega)$ does not. 
}\end{remark}

As in the Euclidean case, if $u$ belongs to $BV_\mathbb{G}(\Omega)$ then its total variation $\|\nabla_{\mathbb G}u\|$ is a finite positive Radon measure and there is a $\|\nabla_{\mathbb G}u\|$-measurable function $\sigma_{u}:\Omega\rightarrow{\mathbb{R}}^{m}$ such that $|\sigma_{u}(x)|=1$ for $\|\nabla_{\mathbb G}u\|$-a.e. $x\in\Omega$ and
\begin{equation}    \label{defsigmaf}
\int_{\Omega}u(x)\mathrm{div}_{\mathbb G}g(x)dx=
\int_{\Omega}\langle g,\sigma_{u}\rangle d\|\nabla_{\mathbb G}u\|
\end{equation}
for all $g\in C_{c}^{1}(\Omega,{\mathbb{R}}^{m})$. 
We denote by $\nabla_{\mathbb G}u$ the vector measure $-\sigma_{u}\|\nabla_{\mathbb G}u\|$, so that $X_{j}u$ is the measure 
$(-\sigma_{u})_{j}\|\nabla_{\mathbb G}u\|$ and
the following integration by parts formula holds true
\begin{equation}   \label{byparts_BV}
\int_{\Omega}u(x)X_{j}g(x)dx = - 
\int_{\Omega}g\left(x\right)d\left( X_{j}u\right) \left( x\right)
\end{equation}
for all $g\in C_c^{1}(\Omega)$. 

\begin{definition}[\sc Sets of finite $\mathbb{G}$-perimeter]\label{perimeter}
If $\chi_{E}$ is the characteristic function of a mea\-su\-rable set $E\subset\mathbb{R}^{N}$, we say that $E$ is a set of finite
$\mathbb{G}$-perimeter in $\Omega$ if $\|\nabla_{\mathbb G}\chi_E\|(\Omega)$ is finite, and we call (generalized inward) 
$\mathbb{G}$-normal the $m$-vector 
\[
\nu_{E}\left(  x\right)  =-\sigma_{\chi_{E}}\left(  x\right)  .
\]
\end{definition}
As customary, we write $P_{\mathbb G}(E)$ instead of $\|\nabla_{\mathbb G}\chi_{E}\|$, $P_{\mathbb G}(E,F)$ instead of $\|\nabla_{\mathbb G}\chi_{E}\|(F)$ for any Borel set $F$. Recall that $\left\vert \nu_{E}\left(  x\right)  \right\vert =1$ for $P_{\mathbb G}(E)$-a.e. $x\in\mathbb{R}^{N}$, so that (\ref{defsigmaf}) takes the form 
\begin{equation}   \label{div thm}
\int_{E}\mathrm{div}_{\mathbb G}g(x)dx=-
\int_{\Omega}\langle g,\nu_{E}\rangle dP_{\mathbb G}(E),
\quad g\in C^1_c(\Omega,\R^m). 
\end{equation}
If $E$ has a smooth boundary, we can compare the generalized ${\mathbb G}$-normal with the Euclidean one, see \cite{CapDanGar94The}, formula (3.2) and \cite[Remark 2.6]{BraMirPal}.
\begin{remark}\label{Esmooth}
{\rm 
If $E$ is a bounded smooth domain in ${\mathbb R}^{N}$ and $n_{E}$ is the Euclidean unit inner normal at $\partial E,$ consider the $m$-vector $v$ whose $j$-th component is defined by
\[
v_j = \sum_{k=1}^{N} \varphi_k^j(x) (n_E)_k(x), \ j=1,\ldots,m,
\]
where the $\varphi_{k}^{j}$ are the coefficients of the vector fields $X_j$ defined in \eqref{defX}. Then
\[
\int_{E}\mathrm{div}_{\mathbb G}g(x)dx=-\int_{\partial E}
\langle g, v \rangle d{\mathcal H}_e^{N-1}\left(  x\right)
\]
(where ${\mathcal H}_e^{N-1}$ denotes the $(N-1)$-dimensional Euclidean Hausdorff measure), from which we read that in this case
\begin{equation}\label{smooth}
\nu_{E} =\frac{v}{\left\vert v \right\vert }, \qquad
P_{\mathbb G}(E) =\left\vert v \right\vert
\left({\mathcal H}_e^{N-1}\mres{\partial E}\right)
\end{equation}
at noncharacteristic points, which are those points of the boundary where $v\neq0$.
}\end{remark}

Let us come to some finer properties of $BV_{\mathbb G}$ functions. 
In order to put formula \eqref{div thm} in a form closer to the classical one, we introduce the notion of {\em measure-theoretic} or {\em essential} boundary, which is a subset of the topological boundary. 
\begin{definition}[\sc Essential boundary]\label{DefEssBdry}
Let $E\subset {\mathbb G}$ be a measurable set. We say that $x\in \partial_{\mathbb G}^*E$ if 
\[
\limsup_{r\to 0}\frac{\lambda_N(E\cap B(x,r))}{\lambda_N(B(x,r))}>0,\qquad
\limsup_{r\to 0}\frac{\lambda_N(B(x,r)\setminus E)}{\lambda_N(B(x,r))}>0 
\]
and we call $\partial_{\mathbb G}^*E$ the {\em measure-theoretic} or {\em essential} boundary of $E$.
\end{definition}

Observe that two different but equivalent distances on ${\mathbb G}$ give the same essential boundary.  

Let us see that the divergence theorem \eqref{div thm} can be rewritten in a form much closer to the classical formula, see \cite[Theorems 5.3, 5.4]{ambr2}, where the problem is settled in general metric measure spaces, see also \cite[Theorem 1.6]{AmbSci10AMPA}.

\begin{theorem}  \label{hausdorffrep}
Given a set of finite perimeter $E\subset{\mathbb G}$, for $P_{\mathbb G}(E,\cdot)$-a.e. $x\in{\mathbb G}$ there is $\bar{r}(x)>0$ such  that 
\[
\ell_{\mathbb G} r^{Q-1} \leq P_{\mathbb G}(E, B_c(x,r)) \leq L_{\mathbb G}r^{Q-1}
\]
for every $r<\bar{r}(x)$, where $0<\ell_{\mathbb G} \leq L_{\mathbb G}<\infty$ are two constants depending only on the group. As a consequence, $P_{\mathbb G}(E,\cdot)$ is concentrated on $\partial_{\mathbb G}^*E$, i.e., $P_{\mathbb G}(E,{\mathbb G}\setminus\partial_{\mathbb G}^*E)=0$. Moreover, 
there is a Borel function 
$\beta_E:{\mathbb R}^N\to [\ell_{\mathbb G},L_{\mathbb G}]$ such that
\begin{equation}\label{intrepr}
P_{\mathbb G}(E,B)=\int_{B\cap \partial_{\mathbb G}^*E}\beta_E (x)\,d{\mathcal S}_{\mathbb G}^{Q-1}(x),\quad \forall B\in {\mathcal B}({\mathbb G}).
\end{equation}
\end{theorem}
The above theorem allows us to rewrite formula \eqref{div thm} as an integral on the essential boundary with respect to the $(Q-1)$-dimensional spherical Hausdorff measure as follows:
\begin{equation}\label{div thm boundary}
\int_{E}\mathrm{div}_{\mathbb G}g(x)dx=
-\int_{\partial_{\mathbb{G}}^{\ast}E}\langle
g,\nu_{E}\rangle \, \beta_E(x)\, d{\mathcal S}_{\mathbb G}^{Q-1} .
\end{equation}

\begin{remark}\label{smoothcase}{\rm We collect here some useful results proved by Franchi, Serapioni and Serra Cassano \cite[Theorem 2.3.5]{fraserser} on functions belonging to $C^1_{\mathbb G}(\Omega)$, for which much more information is available. 

If $\Omega$ is bounded, a function $u$ in $C^1_{\mathbb G}(\Omega)$ also belongs to $BV_{\mathbb G}(\Omega)$ and by \eqref{byparts_BV} the equalities
\[
\int_\Omega X_j^*g(x)u(x)dx = \int_\Omega g(x)X_ju(x)dx, \qquad j=1,\ldots,m,
\]
hold for every $g\in C^1_c(\Omega)$. Recalling \eqref{e-adj}, we find that the measure derivative of $u$ is $\nabla_{\mathbb G}u\,\lambda_N$, where $\lambda_N$ is the Lebesgue measure. Moreover, we say that $S\subset\Omega$ is a ${\mathbb G}$-regular surface if for any $p\in S$ there are an open neighbourhood $U$ of $p$ and $f\in C^1_{\mathbb G}(U)$ such that 
\[
S\cap U = \{x\in U:\ f(x)=0\ \textrm{and}\ \nabla_{\mathbb G}f(x)\neq 0\}.
\]
Let $\Omega$ be an open subset of ${\mathbb R}^N$, $f\in C^1_{\mathbb G}(\Omega)$, $E=\{f<0\}$, $S=\{f=0\}$, and let $p\in\Omega$ be such that $f(p)=0$ and $\nabla_{\mathbb G}f(p)\neq 0$. Then, as proved in \cite[Theorem 2.1]{fraserser2}, there is a neighborhood $U$ of $p$ such that $S\cap U$ has finite perimeter and 
\[
\nu_E(x)=-\frac{\nabla_{\mathbb G}f(x)}{|\nabla_{\mathbb G}f(x)|}, \qquad  x\in S\cap U.
\]
In such a situation the equality $\partial^*_{\mathbb G}(E\cap U) = \partial(E\cap U)$ holds, see \cite[Theorem 3.3]{fraserser2}. Notice also that the topological dimension of a $C^1_{\mathbb G}$-regular surface is $N-1$, see \cite[Proposition 3.1]{fraserser2}, whereas its Hausdorff dimension with respect to the distance $d_\infty$ (or any other equivalent metric) is $Q-1$, see \cite[Corollary 3.7]{fraserser2}. 
}\end{remark}

\begin{remark}\label{thetacostante}
{\rm If $E$ is a finite perimeter set and $\partial_{\mathbb G}^*E$ is ${\mathbb G}$-regular, then formulas \eqref{intrepr} and \eqref{div thm boundary} become simpler. Indeed, in this case for every $x\in \partial_{\mathbb G}^*E$ the normal unit vector $\nu(x)$ is defined and, denoting by $\nu^\bot (x)$ the codimension 1 subspace of $V_1$ orthogonal to $\nu(x)$, we can introduce the hyperplane $N(x)=\nu^\bot (x)\oplus V_2\oplus \cdots \oplus V_\nu$ in $\R^N$. Then, $\beta_E$ is given by
\[
\beta_E(x)=\theta^{-1}
\max_{z\in B(0,1)}\{\H_e^{N-1}(B(z,1)\cap N(x))\}
\]
and thus $\beta_E(x)=1$ for every $x\in\partial^*E$, by Theorem 4.1 in \cite{MagnaniIUMJ} and the definition of the constant $\theta$ in \eqref{deftheta}. This is the reason why we have chosen the distance $d_\infty$. These considerations are important in our proof of Theorem \ref{th-2}, where \eqref{div thm boundary} is applied to sets with finite perimeter such that \emph{a part of} the essential boundary is ${\mathbb G}$-regular. Indeed, Theorem 4.1 in \cite{MagnaniIUMJ} is local, hence if $F\subset\partial_{\mathbb G}^*E$ is ${\mathbb G}$-regular and relatively open, then $\beta_E=1$ in $F$. 
}\end{remark}

We end this subsection with the coarea formula for $BV_{\mathbb G}$ functions, and refer to \cite[Theorem 2.3.5]{fraserser} for its proof. 

\begin{proposition} [\sc Coarea formula in ${\mathbb G}$]\label{coarea} 
If $u\in BV_{\mathbb{G}}(\Omega)$ then for a.e.
$\tau\in\mathbb{R}$ the set $E_{\tau}=\{x\in\Omega:\ u(x)>\tau\}$ has finite $\mathbb{G}$-perimeter and
\begin{equation}    \label{coareaformula}
\left\Vert \nabla_{\mathbb G}u\right\Vert (\Omega)=
\int_{-\infty}^{+\infty}\|\nabla_{\mathbb G}\chi_{E_{\tau}}\|(\Omega)d\tau.
\end{equation}
Conversely, if $u\in L^{1}(\Omega)$ and 
$\int_{-\infty}^{+\infty}\|\nabla_{\mathbb G}\chi_{E_{\tau}}\|(\omega)d\tau<\infty$ 
then $u\in BV_\mathbb{G}(\Omega)$ and equality \eqref{coareaformula} holds. Moreover, if $g:\Omega\to\mathbb{R}$ is a Borel function, then
\begin{equation}    \label{coareag}
\int_{\Omega}g(x)d\left\Vert \nabla_{\mathbb G}u\right\Vert (x)=
\int_{-\infty}^{\infty}\int_{\Omega}g(x)d\|\nabla_{\mathbb G}\chi_{E_{\tau}}\|(x)d\tau. 
\end{equation}
\end{proposition}

\subsection{Proof of Theorem \ref{th-2}}

The proof of Theorem \ref{th-2} is similar to that of Theorem \ref{th-1}. We sketch it and underline the points where different arguments are needed. 

\begin{proof} {\sc of Theorem \ref{th-2}.} Let $\Omega$ be an open subset of $\R^{N}$, and let $u$ be a classical solution to $\E u = f$ in $\Omega$. Let $x_0 \in \Omega$ and let $r_0>0$ be such that $\overline{\Omega_{r_0}(x_0)} \subset \Omega$. Consider an open set $\Omega_1$ such that $\overline{\Omega_{r_0}(x_0)} \subset \Omega_1$ and that $\overline{\Omega_1} \subset \Omega$. Then there exists the fundamental solution $\Gamma^*$ of the equation $\E^* v = 0$ in $\Omega_1$. According to Remark \ref{smoothcase}, $\Gamma^*\in BV_{\mathbb G}(\Omega_1)$, hence we can apply the coarea formula \eqref{coareag} to the set $E_r:= \Omega_r(x_0) =\{\Gamma^*(\cdot,x_0)>r^{2-Q}\}$, which has finite perimeter for a.e. $0<r<r_0$. 

For such a choice of $r$, we choose a positive parameter $\varepsilon>0$ small enough to have $B_\varepsilon(x_0):= \{d_\infty(x,x_0)<\varepsilon\}\subset\Omega_r(x_0)$. We let $\phi_\varepsilon \in C_c^\infty(\R^N)$ be a function such that 
${\rm supp}(\phi)\subset B_\varepsilon(x_0)$ and that $\phi_\varepsilon(x) = 1$ for every $x$ belonging to $B_{\varepsilon/2}(x_0)$. Clearly, $\phi_\varepsilon u \in C_{{\mathbb G},c}^2(\Omega_1)$, then the very definition of $\Gamma^*$ yields
\begin{equation} \label{eq-mvf-g1}
u (x_0) = \phi_\e (x_0) u (x_0) = - \int_\Omega \Gamma^* (x,x_0) \E (\phi_\e u) (x) \, dx.
\end{equation}
We next consider the functions $w(x) := (1 - \phi_\e(x) )u(x)$ and $v(x) := \Gamma^*(x,x_0) - r^{2-Q}$ and we note that 
\begin{equation} \label{eq-div-E*}
\begin{split}
 w(x) \E^* v(x) - v(x) \E w(x) = 
 & \div_{\mathbb G} \big( w(x) A(x) \nabla_{\mathbb G} v(x) 
 - v(x) A(x) \nabla_{\mathbb G} w(x) \big) - 
 \\
 & \div_{\mathbb G} \big( w(x) v(x) b(x) \big),
\end{split}
\end{equation}
which can be written as follows
\begin{equation*}
- v(x) \E w(x) =  \div_{\mathbb G} \, \Phi (x), \qquad 
\Phi(x) := \big( w A \nabla_{\mathbb G} v 
- v A \nabla_{\mathbb G} w - wv b\big)(x).
\end{equation*}
By our choice of $\phi_\e$ we have that $\Phi$ vanishes in $B_{\varepsilon/2}(x_0)$, then it can be extended to a $C_{\mathbb G}^1(\Omega)$ function by setting $\Phi(x_0)= 0$.  We then apply \eqref{div thm boundary} to $\Omega_r(x_0)$ and we find
\begin{equation} \label{eq-div-gk}
 \int_{\Omega_r(x_0)} v(x) \E w(x) dx
=  \int_{\partial_{\mathbb G}^*\Omega_r(x_0)} \scp{\nu,\Phi} \beta_{\Omega_r(x_0)}\, d \S_{\mathbb G}^{Q-1}.
\end{equation}
Here $\nu(x)$ is the generalized unit normal as in Definition \ref{perimeter}, which coincides with 
$\tfrac{\nabla_{\mathbb G} \Gamma^*(x,x_0)}{\abs{\nabla_{\mathbb G} \Gamma^*(x,x_0)}}$ whenever $\nabla_{\mathbb G}\Gamma^*(x,x_0)\neq 0$ according to Remark \ref{smoothcase}. 

Arguing as in the proof of Theorem \ref{th-1}, from the equality $u = w + \phi_\e u$ we deduce 
\begin{equation} \label{eq-div-gk2}
 u (x_0) = - \int_{\Omega_r(x_0)} \left( \tfrac{1}{r^{Q-2}} + v(x) \right) \E (\phi_\e u) (x) dx
\end{equation}
We then recall that $f = \E u = \E w + \E (\phi_\e u)$ so that, by adding \eqref{eq-div-gk2} and \eqref{eq-div-gk}, we find
\begin{equation} \label{eq-div-geps}
 u (x_0) = - \int_{\Omega_r(x_0)} v(x) f(x) + \tfrac{1}{r^{Q-2}}  \E( \phi_\e u)(x) \, dx
 + \int_{\psi_r(x)}\scp{\nu,\Phi} \beta_{\Omega_r(x_0)}\, d\S_{\mathbb G}^{Q-1}.
\end{equation}
Finally, notice that $v=0$ on $\psi_r(x_0)$ and then $\Phi(x)=0$ for every $x\in\partial_{\mathbb G}^*\Omega_r(x_0)$ such that 
$\nabla_{\mathbb G}\Gamma^*(x,x_0)=0$, hence we can write
\begin{align*}
\int_{\partial_{\mathbb G}^*\Omega_r(x_0)} &
\scp{\nu,\Phi} \beta_{\Omega_r(x_0)}\, d \S_{\mathbb G}^{Q-1} = 
\int_{\psi_r(x_0)\setminus\{\nabla_{\mathbb G}\Gamma^*(\cdot,x_0)=0\}} \scp{\nu,\Phi}\, d \S_{\mathbb G}^{Q-1}
\\
&=
\int_{\psi_r(x_0)\setminus\{\nabla_{\mathbb G}\Gamma^*(\cdot,x_0)=0\}} \Big\langle\frac{\nabla_{\mathbb G}\Gamma^*(x,x_0)}{|\nabla_{\mathbb G}\Gamma^*(x,x_0)|},A\nabla_{\mathbb G}\Gamma^*(x,x_0)\Big\rangle u(x) \, d \S_{\mathbb G}^{Q-1}
\\
&=
\int_{\psi_r(x_0)} K_{\mathbb G}(x_0,x)u(x) \, 
d \S_{\mathbb G}^{Q-1},
\end{align*}
where we have taken into account that $\beta_{\Omega_r(x_0)}=1$ in the ${\mathbb G}$-regular part of $\psi_r(x_0)$, see Remark \ref{thetacostante}. The last step to conclude the proof of \eqref{e-meanvalue-g} consists in letting $\varepsilon \to 0$ as in the proof of Theorem \ref{th-1} to get rid of the integral of $r^{2-Q}\E(\phi_\e u)$ in \eqref{eq-div-geps}.

To deduce \eqref{th-2-2} from \eqref{e-meanvalue-g} we argue as in Section \ref{SecUnif}, replacing $N$ by $Q$ and using the coarea formula provided by Theorem \ref{coarea} in the last step. 
\end{proof}

\begin{remark}\label{regularGamma}
{\rm If the coefficents $a_{ij}$ are $C^\infty$ then the fundamental solution $\Gamma^*$ is $C^\infty$ as well and the level sets $\{\Gamma^*>c\}$ are smooth surfaces for almost all $c\in\R$. Therefore, we may take into account Remark \ref{Esmooth} and write the surface integral in \eqref{e-meanvalue-g} by the simpler form
\[
\int_{\psi_r(x_0)} K_{\mathbb G} (x_0,x) u(x) \, d{\mathcal H}_e^{N-1}(x),
\]
i.e., we may use the $(N-1)$-dimensional Euclidean Hausdorff measure.
}\end{remark}

\noindent {\bf Acknowledgments.} This research was partially supported by the grant of Gruppo Naziona\-le per l’Analisi Matematica, la Probabilità e le loro Applicazioni (GNAMPA) of the Istituto Nazionale di Alta Matematica (INdAM). We thank Nicola Garofalo and Ermanno Lanconelli for their interest in our work and Valentino Magnani, Francesco Serra Cassano and Davide Vittone for several useful discussions.


\begin{thebibliography}{10}

\bibitem{ambr2}
{\sc L. Ambrosio}, {\em Fine properties of sets of finite perimeter in doubling metric measure spaces}, Set-Valued Anal., 10, (2002), 111-128.

\bibitem{AmbSci10AMPA}
{\sc L. Ambrosio, M. Scienza}, 
{\em Locality of the perimeter in Carnot groups and chain rule},
Ann. Mat. Pura Appl. 189 (2010), 661-678.

\bibitem{BLU2}
{\sc A.~Bonfiglioli, E.~Lanconelli, F.~Uguzzoni}, {\em Fundamental solutions for non-divergence form operators on stratified groups}, 
Trans. Amer. Math. Soc., 356 (2003), 2709-2737.

\bibitem{LibroBLU}
{\sc A.~Bonfiglioli, E.~Lanconelli, F.~Uguzzoni}, {\em Stratified {L}ie {G}roups and {P}otential {T}heory for their sub--{L}aplacians}, Springer Monographs in Mathematics 2007.

\bibitem{BraMirPal}
{\sc M. Bramanti, M. Miranda, D. Pallara}, {\em Two characterization of $BV$ functions on Carnot groups via the heat semigroup}, Int. Math. Res. Not., 17 (2012), 3845-3876.

\bibitem{CapDanGar94The}
{\sc L. Capogna, D. Danielli, N.Garofalo}, {\em The geometric Sobolev embedding for vector fields and the isoperimetric inequality}, 
Comm. Anal. Geom., 2 (1994), 203-215.

\bibitem{CittiGarofaloLanconelli}
{\sc G.~Citti, N.~Garofalo, E.~Lanconelli}, {\em Harnack's inequality for sum of squares of vector fields plus a potential}, Amer. J. Math. 115 (1993), 699-734.

\bibitem{FabesGarofalo}
{\sc E.~B. Fabes, N.~Garofalo}, {\em Mean value properties of solutions to parabolic equations with variable coefficients}, 
J. Math. Anal. Appl., 121 (1987), 305-316.

\bibitem{Folland}
{\sc G.~B. Folland}, {\em Subelliptic estimates and function spaces on nilpotent {L}ie groups}, Ark. Mat., 13 (1975), pp.~161--207.

\bibitem {fraserser}
{\sc B. Franchi, R. Serapioni, F. Serra Cassano}, 
{\em Meyers-Serrin type theorems and relaxation of variational integrals depending on vector fields}, 
Houston J. Math. 22 (1996), 859-890.

\bibitem {fraserser2}
{\sc B. Franchi, R. Serapioni, F. Serra Cassano}, 
{\em Regular hypersurfaces, intrinsic perimeter and implicit function theorem in Carnot groups},
Comm. Anal. Geom., 11 (2003), 909-944.

\bibitem {fraserser3}
{\sc B. Franchi, R. Serapioni, F. Serra Cassano}, 
{\em On the strcucture of finite perimeter sets in Carnot step 2 groups},
J. Geom. Anal., 13 (2003), 431-466.

\bibitem{GarofaloLanconelli-1989}
{\sc N.~Garofalo, E.~Lanconelli}, {\em Asymptotic behavior of fundamental solutions and potential theory of parabolic operators with variable coefficients}, Math. Ann., 283 (1989), 211--239.

\bibitem{GarofaloLanconelli-1990}
{\sc N.~Garofalo, E.~Lanconelli}, {\em Level sets of the fundamental solution and Harnack inequality for degenerate equations of Kolmogorov type}, Trans. Amer. Math. Soc. 321 (1990), 775-792.

\bibitem{Kalf}
{\sc H. Kalf}, {\em On E. E. Levi's method of constructing a fundamental solution for second-order elliptic equations}, Rend. Circ. Mat. Palermo, (2) 41 (1992), 251-294. 

\bibitem{KogojLanconelli2}
{\sc A.~E. Kogoj, E.~Lanconelli}, {\em An invariant {H}arnack inequality for a class of hypoelliptic ultraparabolic equations}, Mediterr. J. Math., 1 (2004), pp.~51--80.

\bibitem{MagnaniIUMJ}
{\sc V. Magnani}, {\em A New differentiation, shape of the unit ball, and perimeter measure}, Indiana Univ. Math. J., 66 (2017), pp. 183-204.

\bibitem{MagnaniCalcVar}
{\sc V. Magnani}, {\em Towards a theory of area in homogeneous groups}, Calc. Var. Partial Diff. Eq. 58 (2019), Paper No. 91, 39 pp.

\bibitem{MalPalPol}
{\sc E.~Malagoli, D.~Pallara, S.~Polidoro}, {\em Mean value formulas for classical solutions to uniformly parabolic equations in divergence form}, to appear. 

\bibitem{Morgan}
{\sc F. Morgan}, {\em Geometric measure theory. A beginner's guide, $3^{rd}$ ed.}, Academic Press, 2000.

\bibitem{PP}
{\sc D.~Pallara, S.~Polidoro}, {\em Mean value formulas for classical solutions to subelliptic equations in stratified Lie groups}, to appear. 

\bibitem{Pini1951}
{\sc B.~Pini}, {\em Sulle equazioni a derivate parziali, lineari del secondo ordine in due variabili, di tipo parabolico}, Ann. Mat. Pura Appl. (4), 32 (1951), pp.~179--204.

\bibitem{Polidoro2}
{\sc S.~Polidoro}, {\em On a class of ultraparabolic operators of {K}olmogorov-{F}okker-{P}lanck type},
Le Matematiche (1), 49 (1994), 53-105.

\bibitem{SerCassEMS}
{\sc F. Serra Cassano}, {\em Some topics of geometric measure theory in Carnot groups}, in: Geometry, Analysis and Dynamics on sub-Riemannian
Manifolds, Volume I, D. Barilari, U. Boscain, M. Sigalotti eds., EMS 2016, 1-121.

\bibitem{Varadhan1967behavior}
{\sc S.~R.~S. Varadhan}, {\em On the behavior of the fundamental solution of the heat equation with variable coefficients}, Communications on Pure and Applied Mathematics 20~(2) (1967), 431-455.

\bibitem{Varadhan1967diffusion}
{\sc S.~R. Varadhan}, {\em Diffusion processes in a small time interval}, Commu. Pure and Appl. Math. 20~(4) (1967), 659-685.

\bibitem{Watson1973}
{\sc N.~A. Watson}, {\em A theory of subtemperatures in several variables}, Proc. London Math. Soc. (3), 26 (1973), pp.~385--417.
  
\end{thebibliography}
\end{document}